# Quadratically-Regularized Optimal Transport on Graphs

Montacer Essid, Justin Solomon

**Abstract**

Optimal transportation provides a means of lifting distances between points on a geometric domain to distances between signals over the domain, expressed as probability distributions. On a graph, transportation problems can be used to express challenging tasks involving matching supply to demand with minimal shipment expense; in discrete language, these become minimum-cost network flow problems. Regularization typically is needed to ensure uniqueness for the linear ground distance case and to improve optimization convergence; state-of-the-art techniques employ entropic regularization on the transportation matrix. In this paper, we explore a quadratic alternative to entropic regularization for transport over a graph. We theoretically analyze the behavior of quadratically-regularized graph transport, characterizing how regularization affects the structure of flows in the regime of small but nonzero regularization. We further exploit elegant second-order structure in the dual of this problem to derive an easily-implemented Newton-type optimization algorithm.

## 1 Introduction

Since its formalization in the eighteenth century, the *optimal transportation* problem has provided practical and theoretical challenges in a variety of scientific, mathematical, and engineering disciplines. Based on the natural problem of optimizing a matching between supply and demand to minimize shipping costs, optimal transportation has been proposed in countless fields under different names: the Monge–Kantorovich problem in economics [42], the Hitchcock–Koopmans problem in engineering optimization [18], the earth mover's distance in computer vision [31], the Wasserstein distance in analysis [42], minimum-cost network flow in graph theory [19], and the Mallows distance in statistics [27]—to name a few. This appearance and reappearance of optimal transportation underscores its fundamental nature and value as a tool for modeling and analysis.

Modern computational applications of transport impose extreme demands on scalability. Whereas problems involving a relatively small number of sources and targets can be solved using standard linear programming, extending to thousands or millions of nodes pushes the limits of generic polynomial-time algorithms. An additional challenge is provided by nonuniqueness of the solution to transport with linear ground distance, leading to unpredictable output by standard software packages; this is exacerbated on large graphs in which cycles and other features make multiple shortest paths and matchings indistinguishable. Hence, recent incarnations of computational transport must identify additional structure in the problem that can be leveraged to overcome issues in the most generic instances of the problem.

In this paper, we consider a critical application of "structured" optimal transportation, namely when the transport cost comes from distances along an underlying graph. This realistic assumption appears in many contexts, e.g. transporting goods along a road map or matching servers to clients on a computer network. While classical network flow-style phrasings of this problem have been tackled in the algorithms community, here we take inspiration from numerical methods approximating transport over continuum domains to study *regularized* transport over a network; regularization ensures a unique solution to the transport and smooths the objective landscape of the problem. We provide theoretical characterization of the behavior of the regularized problem as it deviates from the non-regularized case. Furthermore, the appearance of the graph Laplacian in the dual problem allows us to propose a Newton-style optimization algorithm for the problem that leverages sparsity and low-rank structure to invert the Hessian efficiently.

## 2 Related Work

Transport over graphs has a rich history and has been considered in mathematics, operations research, computer science, and many other disciplines. At the broadest level, this problem is a slight generalization of the *linear assignment* problem and can be solved using several classical techniques; see e.g. [6, 4] for discussion.

Without regularization, the particular problem we study would be an instance of *minimum-cost flow without edge capacities*. [30] discusses classical algorithms for this linear programming problem, such as the cycle canceling [20], network simplex [28], and Ford–Fulkerson algorithms [15]. Classical methods such as these often are not accompanied with systematic "tie-breaking" strategies in case the network flow is non-unique, indicating the potential application of a strictly convex variation of the problem such as ours.

The theoretical computer science community recently has reconsidered this class of problems from an optimization perspective. In the undirected, capacity-free case, [36] proposes a preconditioner for approximate minimum-cost flow problems that achieves nearly-linear runtimes. [10] considers the more general case of minimum-cost flow in directed graphs with unit capacities, extending a framework proposed in [26] for approaching graph-based problems using the interior point method. These algorithms are primarily of theoretical interest but do employ interior point-style methods, possibly implying a systematic choice of flows in the case of multiple optima.

In the continuum, the theory of optimal transport [42] classifies problems structured similarly to minimum-cost flow as the *1-Wasserstein distance* or *Beckmann* problem [1]. See [33] for analysis and [42, 17, 32] for theoretical discussion. Even over general spaces, solutions of the 1-Wasserstein problem generally are nonunique and include some degenerate optima [42, §2.4.6]. Numerical algorithms for these problems include [39], which uses finite element methods for a vector field version of this problem accelerated using spectral decompositions, and [21], which adapts primal–dual methods designed for $L^1$-regularized optimization.

Our work involves a *regularized* model of transport, in which the cost of a per-edge flow is augmented with a strictly convex term. For bipartite graphs, by far the most popular regularized approach to transport involves entropic regularization [11, 3], which leads to an instance of the well-known Sinkhorn–Knopp rescaling algorithm [38], also known as the iterative proportional fitting procedure [44, 13]. These algorithms are extremely effective for bipartite transport problems thanks to the elegant algorithms in this case, but to our knowledge no Sinkhorn-like method has been formulated for transport over more general graphs.

We instead apply quadratic regularization, accompanying the $L^1$-style transport objective function with an additional $L^2$ term. This regularizer was used in [21, 22, 23] to derive a parallelizable primal-dual algorithm for regularized 1-Wasserstein distance, on regular grids of $\mathbb{R}^d$. Our work is complementary to theirs, providing algorithm-independent analysis of the structure of flow in the presence of quadratic regularization, in particular how its coefficient controls sparsity, as well as analysis of the dual yielding a Newton-type optimization method; the discussion here provides fine-grained information about the quadratic case, whereas many of their constructions apply to other regularizers.

[43] studies the general case of $L^1$ minimization problems with an $L^2$ regularization term. In particular, they establish a similar result on sparsity of the regularized solutions for a broader class of problems than the ones studied in this paper; our result can be obtained by applying their results to a graph divergence operator. In the context of the minimum cost flow problem, our paper obtains the result by taking a different proof approach: We analyze changes in mass flows in the primal solution

when the regularization parameter changes instead of studying the dual problem.

In the bipartite case, we could replace entropic regularization with quadratic regularization for a less efficient but Sinkhorn-like algorithm through the alternating projection framework of [3]; these projections would require sorting a list of floating-point values in each iteration [14].

While our theoretical and practical consideration is largely focused on the small-regularization regime, in the case of high regularization our flows begin to resemble *electrical flows* on graphs [8, 26]. See e.g. [9] for relevant notions from spectral graph theory.

Finally, we note a few recent works [25, 40] with alternative models for optimal transport inspired by fluid dynamics interpretations of transport with quadratic costs [2]. These methods hold some potential to overcome nonuniqueness issues associated with linear transport costs but relatively few faithful numerical discretization and optimization techniques exist; [16] suggests one possible approach. One recent paper applies regularization to the dynamical problem [24], potentially suggesting an alternative means to introduce regularization in optimal transport.

## 3 Quadratically-Regularized Transport on Graphs

### 3.1 Graph Transport without Regularization

Suppose $G = (V, E)$ is a connected graph with directed edges $E \subseteq V \times V$. We associate edges in $E$ with a vector of edge weights $c \in \mathbb{R}_+^{|E|}$.

Denote $\mathrm{Prob}(V)$ to be the probability simplex over $V$, that is

$$\mathrm{Prob}(V) := \{\rho \in \mathbb{R}^{|V|} : \mathbb{1}^\top \rho = 1 \text{ and } \rho \geq 0\}.$$

The *1-Wasserstein* [42], *optimal transportation* [42], or *earth mover's* [31] distance between two distributions $\rho_0, \rho_1 \in \mathrm{Prob}(V)$ is defined as

$$\mathcal{W}_1(\rho_0, \rho_1) := \begin{cases} \min_{T \in \mathbb{R}^{|V| \times |V|}} & \mathrm{Tr}(C^\top T) \\ \text{s.t.} & T \geq 0 \\ & T\mathbb{1} = \rho_0 \\ & T^\top \mathbb{1} = \rho_1. \end{cases} \tag{1}$$

Here, the unknown matrix $T \in \mathbb{R}^{|V| \times |V|}$ is a *transportation plan* transforming $\rho_0$ into $\rho_1$. The precomputed cost matrix $C \in \mathbb{R}_+^{|V| \times |V|}$ contains shortest-path distances between each pair of vertices on the graph given edge lengths in $c$. Intuitively, over all possible transportation matrices $T$, we minimize the total "work" measured as the sum of mass $T_{vw}$ moved from vertex $v$ to vertex $w$ times distance $C_{vw}$. Note that $C$ needs not be symmetric since the edges in $E$ are directed.

The cost of moving mass between two vertices $i$ and $j$ on $G$ can be decomposed as the cost of moving mass along each edge on the shortest path between $i$ and $j$. Formalizing this argument provides an alternative to (1) with one variable per edge in $E$, a considerable savings when $G$ is sparse. Define an incidence matrix $D \in \{-1, 0, 1\}^{|E| \times |V|}$ as

$$D_{ev} := \begin{cases} -1 & \text{if } e = (v, w) \text{ for some } w \in V \\ 1 & \text{if } e = (w, v) \text{ for some } w \in V \\ 0 & \text{otherwise.} \end{cases} \tag{2}$$

This matrix is an analog of the gradient operator for functions on $\mathbb{R}^n$. Then, an alternative formula for $\mathcal{W}_1$ is

$$\mathcal{W}_1(\rho_0, \rho_1) := \begin{cases} \min_{J \in \mathbb{R}^{|E|}} & \sum_{e \in E} c_e J_e \\ \text{s.t.} & J \geq 0 \\ & D^\top J = \rho_1 - \rho_0. \end{cases} \tag{3}$$

The vector $J$ contains one directed flow per edge; the constraint expresses the requirement that $J$ flows from $\rho_0$ to $\rho_1$. The transpose $D^\top$ takes the place of (negative) divergence for vector fields on $\mathbb{R}^n$. In the language of smooth optimal transport, this transformation of the $\mathcal{W}_1$ problem for graphs becomes the *Beckmann problem* [34] for measures on $\mathbb{R}^n$. [39] provides optimization techniques and applications of this objective to computer graphics; [21] presents an algorithm targeted to image domains.

## 3.2 Regularized Transport

Regardless of whether we write the problem as (1) or (3), the transport problem as formulated above suffers from nonuniqueness of the optimized variable $T$ or $J$, that is, the problem is convex but not strictly convex. This makes the output of transport algorithms with linear costs unpredictable at best.

A typical approach to making the problem better posed is to add a regularizer, i.e. a second objective term adding stronger convexity to the problem. For instance, [11, 3] regularize the matrix $T$ in (1) by subtracting a term proportional to its entropy $-\sum_{vw \in V} T_{vw} \ln T_{vw}$. This popular work arguably has revitalized interest in computational optimal transportation by providing an efficient approximation to a wide variety of challenging problems in this domain, but entropic regularization is worth reconsidering in the context of the graph-based formulation (3) for a few reasons:

- The regularization is written in terms of elements of $T$ rather than elements of $J$. Entropic regularization on $J$ rather than $T$ is possible but does not appear to admit as elegant an optimization algorithm.

- As the coefficient of the regularizer is decreased to zero, algorithms based on alternating projection exhibit slower convergence as well as numerical issues dealing with near-zero values. Note [35] recently introduces some improvements that help entropy-regularized transport in the small regularization regime.

- No matter how small the entropic regularizer, the resulting transportation matrix $T$ satisfies $T_{vw} > 0$ strictly for *all* pairs of vertices $v, w \in V$. This contrasts qualitatively with sparsity properties for non-regularized transport.

We emphasize the first point above. Entropic regularization appears to be most effective when applied directly to the transport matrix but is less relevant for problems where the variables can be decomposed into per-edge flows.

In our paper, we explore algorithms for an alternative quadratically-regularized model of regularized transport:

$$\mathcal{W}_{1,\alpha}(\rho_0, \rho_1) := \begin{cases} \min_{J \in \mathbb{R}^{|E|}} & \sum_{e \in E} c_e J_e + \frac{\alpha}{2} \sum_e J_e^2 \\ \text{s.t.} & J \geq 0 \\ & D^\top J = \rho_1 - \rho_0. \end{cases} \quad (4)$$

Our reasons for studying (4) include the following:

- The quadratic regularizer allows for $J_e = 0$ exactly. Furthermore, this regularizer is amenable to the analysis in §4 showing that $\alpha$ modulates sparsity of $J$ in a controlled fashion.

- The algorithm we will propose in §5 has little in common with alternating projection techniques used for entropically-regularized transport and is suited to low regularization and sparse graphs with $|E| \ll |V|^2$.

## 3.3 Dual Problem

The optimization problem (4) is a convex quadratic program with affine constraints, and hence it exhibits strong duality by the affine Slater condition [5, eq. (5.27)]. Denote the positive part of a vector $v \in \mathbb{R}^n$ as $v_+$ with elements $(v_+)_i := \max\{v_i, 0\}$, and let $|v|$ denote its Euclidean norm. Then, strong duality implies the following proposition:

**Proposition 1** (Duality). *Quadratically-regularized transport on graphs* (4) *can be computed as follows:*

$$\mathcal{W}_{1,\alpha}(\rho_0, \rho_1) = \frac{1}{\alpha} \sup_{p \in \mathbb{R}^{|V|}} \left[ \alpha f^\top p - \frac{1}{2} |(Dp - c)_+|_2^2 \right], \tag{5}$$

*where* $f := \rho_1 - \rho_0$ *and* $c \in \mathbb{R}^{|E| \times 1}$ *is the vector of costs per edge. Furthermore, the primal variable* $J_e$ *on edge* $e \in E$ *is zero whenever* $(Dp - c)_e \leq 0$.

*Proof.* With the notation defined above,

$$\mathcal{W}_{1,\alpha}(\rho_0, \rho_1) := \begin{cases} \min_{J \in \mathbb{R}^{|E|}} & \sum_{e \in E} c_e J_e + \frac{\alpha}{2} \sum_e J_e^2 \\ \text{s.t.} & J \geq 0 \\ & D^\top J = f \end{cases}$$

$$= \min_{J \in \mathbb{R}^{|E|}_+} \max_{p \in \mathbb{R}^{|V|}} \left[ c^\top J + \frac{\alpha}{2} J^\top J + (f - (D^\top J))^\top p \right]$$

$$= \max_{p \in \mathbb{R}^{|V|}} \left[ f^\top p + \min_{J \in \mathbb{R}^{|E|}_+} \left( J^\top (c - Dp) + \frac{\alpha}{2} J^\top J \right) \right]. \tag{6}$$

Here, we switched max and min by strong duality, as noted above.

Now, the inner minimum can be explicitly computed since it is a quadratic function of the variable $J$; optimality yields the complimentary slackness condition

$$J = \frac{(Dp - c)_+}{\alpha}.$$

Substituting into (6) provides the desired unconstrained dual formulation:

$$\mathcal{W}_{1,\alpha}(\rho_0, \rho_1) = \sup_{p \in \mathbb{R}^{|V|}} \left[ f^\top p - \frac{1}{2\alpha} |(Dp - c)_+|_2^2 \right]$$

□

This dual problem suggests posing our problem in terms of the *active set* of edges implied by a dual variable $p \in \mathbb{R}^{|V|}$:

**Definition 1** (Active set). *The* active set *of edges in* $E$ *associated to a dual variable* $p$ *is the set* $S(p) := \{e \in E : (Dp - c)_e > 0\}$. *In a minor abuse of notation, we will denote the active set of edges in* $E$ *associated to a primal variable* $J$ *as the set* $S(J) := \{e \in E : J_e > 0\}$.

Given $\alpha > 0$, the solution $J^\alpha$ of the primal problem (4), and the solution $p^\alpha$ of the dual problem (5), by complimentary slackness we know $S(J^\alpha) = S(p^\alpha)$. Hence, we denote $S(\alpha) := S(J^\alpha) = S(p^\alpha)$ to be the active set of edges associated to the optimization problem with regularizing coefficient $\alpha$.

For development of our optimization algorithm in §5, we make a key observation regarding the derivatives of the dual problem (5). For convenience, define $M(p)$ to be a diagonal matrix associated with the active set $S(p)$:

$$M(p)_{ee} := \begin{cases} 0 & \text{if } e \notin S(p) \\ 1 & \text{if } e \in S(p). \end{cases} \tag{7}$$

This allows us to write the objective in (5) as

$$\sup_{p \in \mathbb{R}^{|V|}} \underbrace{\left[ \alpha f^\top p - \frac{1}{2} (Dp - c)^\top M(p)(Dp - c) \right]}_{g(p)}.$$

Differentiating this expression where $M(\cdot)$ is constant gives the expressions

$$\nabla g(p) = \alpha f - D^\top M(p)(Dp - c) \tag{8}$$

$$\text{Hess}[g](p) = -D^\top M(p) D \tag{9}$$

In the language of spectral graph theory [9], the unweighted Laplacian matrix of a graph $G = (V, E)$ is the matrix $L := D^\top D$. The second expression immediately provides a key observation:

**Proposition 2.** *The Hessian of the dual problem* (5) *is the unweighted Laplacian $L(p)$ of the active subgraph $G(p) := (V, S(p))$.*

Note that $L(p)$ is **never full-rank**; in particular, the null space of $L(p)$ is spanned by the indicator vectors of the connected components of $G$. This prohibits direct application of Newton's method to optimizing the dual problem.

## 4 Effect of Regularizer on Active Graph

Transport with quadratic regularization is less well-understood than its entropically-regularized counterpart. A key feature is that this problem allows elements of $J$ to take zero values, while entropically-regularized transport would force $J > 0$. To better understand this and other properties, in this section we provide theoretical characterization of the computed flow $J$ as a function of the regularizing coefficient $\alpha$.

### 4.1 Notation and Basic Properties

Recall that the flow $J$ is nonzero only on the active edges $S(p)$. Denote $p(\alpha)$ to be an optimal dual variable as a function of the regularizing coefficient $\alpha$. Our goal is to provide some intuition for the effect of $\alpha$ on $J$ and $p$.

To streamline discussion, let (LP) denote the non-regularized problem

$$\begin{cases} \min_{J \in \mathbb{R}^{|E|}} & \sum_{e \in E} c_e J_e \\ \text{s.t.} & J \geq 0 \\ & D^\top J = f, \end{cases} \tag{LP}$$

and let $J_0$ be a solution of (LP) with active set $S(J_0)$; note that $J_0$ is nonunique. Furthermore, define the union of all active sets of the non-regularized problem as

$$S_0 = \bigcup_{J_0 \text{ solution of (LP)}} S(J_0). \tag{10}$$

Let (QP) denote the quadratically-regularized primal min-cost flow problem

$$\begin{cases} \min_{J \in \mathbb{R}^{|E|}} & \sum_{e \in E} c_e J_e + \frac{\alpha}{2} J_e^2 \\ \text{s.t.} & J \geq 0 \\ & D^\top J = f, \end{cases} \tag{QP}$$

and let $J^\alpha$ be the solution of (QP) with active set $S(\alpha)$. Using $|J^\alpha|$ to denote the Euclidean norm of vector $J^\alpha$, we can rewrite the objective of (QP) as

$$V_\alpha(J) = c^\top J + \frac{\alpha}{2} |J|^2.$$

$S(J_0)$ might not be unique, as linear programming might exhibit multiple solutions $J_0$ with different active sets. The quadratically regularized problem, however, must exhibit uniqueness of the solution thanks to strict convexity. By standard Γ-convergence arguments, $J^\alpha$ converges to a solution of (LP) as $\alpha \to 0$.

The objective function of the quadratically regularized problem (QP) is the sum of an $L^1$ cost $c^\top J$ and an $L^2$ regularizer $|J|^2$, where the parameter $\alpha$ tunes the influence of the two terms. Intuitively, as $\alpha$ decreases, solutions of (QP) will favor low $L^1$ cost at the expense of the $L^2$ cost, which will increase instead. This behavior is formalized as follows:

**Proposition 3.** *Let $0 < \alpha < \alpha'$, and $J^\alpha$ be the solution of (QP). If $J^\alpha \neq J^{\alpha'}$, then $c^\top J^\alpha < c^\top J^{\alpha'}$ and $|J^\alpha|^2 > |J^{\alpha'}|^2$.*

*Proof.* Since $J^\alpha$ is the unique minimizer of the strictly convex optimization problem, we have $[\nabla V_\alpha(J^{\alpha'})]^\top(J^\alpha - J^{\alpha'}) < 0$, or equivalently

$$c^\top(J^\alpha - J^{\alpha'}) + \alpha(J^{\alpha'})^\top(J^\alpha - J^{\alpha'}) < 0. \tag{11}$$

Reversing the relationship for $\alpha'$ shows

$$c^\top(J^\alpha - J^{\alpha'}) + \alpha'(J^{\alpha'})^\top(J^\alpha - J^{\alpha'}) \geq 0.$$

Subtracting these two expressions yields the inequality

$$(J^{\alpha'})^\top(J^\alpha - J^{\alpha'}) > 0. \tag{12}$$

Substituting this expression into (11) shows our first claim $c^\top(J^\alpha - J^{\alpha'}) < 0$. For our next claim, we simplify

$$|J^\alpha|^2 = |J^{\alpha'} + (J^\alpha - J^{\alpha'})|^2 = |J^{\alpha'}|^2 + 2(J^{\alpha'})^\top(J^\alpha - J^{\alpha'}) + |J^\alpha - J^{\alpha'}|^2$$
$$> |J^{\alpha'}|^2 + |J^\alpha - J^{\alpha'}|^2 \text{ by } (12)$$
$$> |J^{\alpha'}|^2, \text{ as needed.} \qquad \square$$

$\square$

## 4.2 Overview of Main Results

Careful analysis of our problem gives us a more precise description of the evolution of the solution as we decrease $\alpha$. We provide two related results regarding the behavior of $J^\alpha$ for small $\alpha > 0$:

**Proposition 4** (Sparsity). *There exists a constant $\tilde{\alpha} > 0$ depending on the graph $G$ and data $f$ such that for all $\alpha \in (0, \tilde{\alpha})$, the solution $J^\alpha$ of (QP) is also a solution of (LP).*

Resulting from the proposition above is the following corollary:

**Corolary 1** (Sparsity 2). *There exists $\tilde{\alpha} > 0$ depending on the graph $G$ and data $f$ as well as a unique solution $J_0$ of (LP) such that for all $\alpha \in (0, \tilde{\alpha})$, $J_0$ is the unique solution of (QP).*

The following counterexamples rule out stronger theoretical results that have frequently been observed experimentally:

- In Proposition 4, we do not necessarily have all solutions of (LP) are solutions of (QP) for low values of $\alpha$. In particular, it is possible that the solution $J^\alpha$ is *sparser* than *some* solution of (LP), as shown in the counterexample of Figure 1.

- We do not necessarily have monotonicity of the active sets in $\alpha$, that is $S(\alpha) \subseteq S(\alpha')$ for $0 < \alpha < \alpha'$, as shown in the counterexample of Figure 2. A monotonicity conjecture will be formulated in a subsequent section, after introducing appropriate notation.

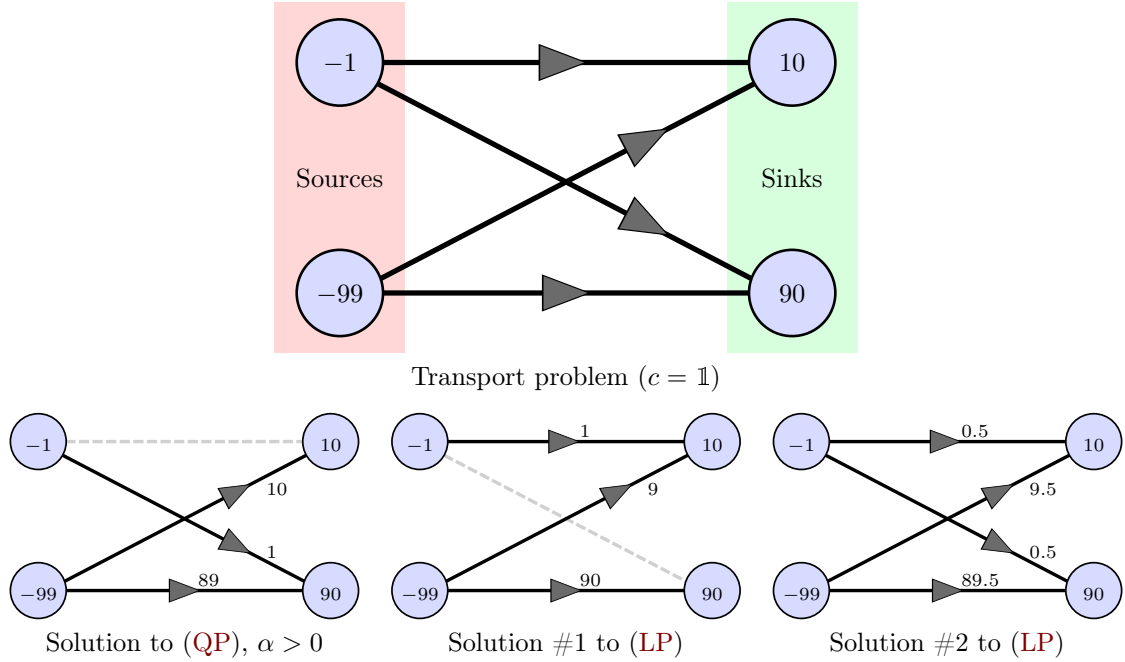

Figure 1: Counterexample showing that the regularized $J^\alpha$ ($\alpha > 0$) can be sparser than *some* nonunique non-regularized solution. Here we show a simple flow problem on a directed graph (top), where values of $f$ are shown on the nodes. When $\alpha > 0$, the solution to (QP) with $\alpha > 0$ has one sparse edge (bottom left), but when $\alpha = 0$ multiple solutions exist with different sparsity patterns (bottom middle/right).

### 4.3 Proofs

#### 4.3.1 Preliminaries

For completeness, we recall standard notions regarding flow decomposition from [30]. We begin with a few definitions:

- A *directed path* is an ordered set of *distinct* nodes $\{v_{i_1} - v_{i_2} - \cdots - v_{i_n}\}$ such that $(v_{i_k}, v_{i_{k+1}}) \in E$ for all $k \in \{0, \ldots, n-1\}$. $\mathcal{P}$ denotes the set of directed paths in $G$.

- For $r \in \mathcal{P}$, define $s(r) := v_{i_1}$ to be the starting node of $r$, $t(r) := v_{i_n}$ to be its ending node, and $\text{len}(r) = n - 1$ to be its length. Define $P := \{(v_{i_1}, v_{i_{k+1}}) : k \in \{1, \ldots, \text{len}(r)\}\} \subseteq E$ as the set of $r$'s edges.

- A *cycle* is an ordered set of nodes $\{v_{i_1} - v_{i_2} - \cdots - v_{i_n}\}$ such that $(v_{i_k}, v_{i_{k+1}}) \in E$ for all $k \in \{0, \ldots, n-1\}$. Furthermore, all nodes are distinct except the starting and ending nodes, which are the same ($v_{i_1} = v_{i_n}$). Let $\mathcal{C}$ denote the set of cycles.

- Given $r \in \mathcal{P} \cup \mathcal{C}$, define a flow $\delta(r) : E \to \{0, 1\}$ via

$$[\delta(r)]_e := \begin{cases} 1 \text{ if } e \in P \\ 0 \text{ if } e \notin P. \end{cases}$$

We will use $\delta(P)$ to refer to the analogous flow given an edge set $P$.

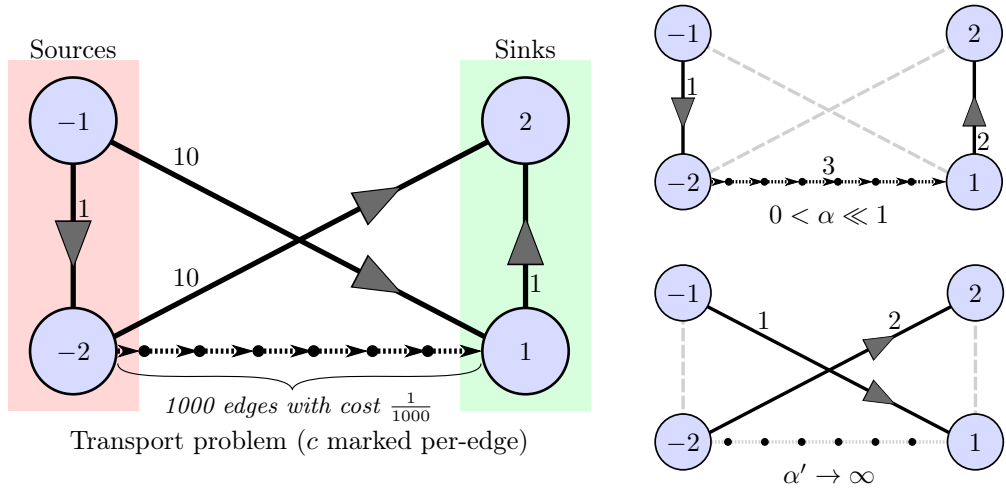

Figure 2: Counterexample showing that the regularized $J^{\alpha'}$ does not necessarily contain the active network of the regularized $J^\alpha$. Here we show a simple flow problem on a directed graph (left), where values of $f$ are shown on the nodes and the cost vector $c$ is shown per edge. When $\alpha$ is close to 0 (top right), the solution to (QP) does not use the high $L^1$ cost edges, but rather transports the mass using the low $L^1$ cost edges. When $\alpha'$ is large (bottom right), it is not advantageous to use the low $L^1$ cost path because of the multiple successive edges that have high $L^2$ cost. The solution prefers then to use the high $L^1$ cost edges.

- Given a flow on paths and cycles $\hat{J} : \mathcal{P} \cup \mathcal{C} \to \mathbb{R}$, we can define an *arc flow* $J : E \to \mathbb{R}$ by summation. For $e \in E$, this flow can be written

$$J_e := \sum_{r \in \mathcal{P}} \delta_e(r) \hat{J}(r) + \sum_{c \in \mathcal{C}} \delta_e(c) \hat{J}(c). \tag{13}$$

We will only need to consider flows $\hat{J}$ satisfying the constraint $D^\top J = f$. The following theorem decomposes every arc flow satisfying the divergence constraint $D^\top J = f$ into path and cycle flows:

**Theorem 1** (Flow decomposition theorem [30, Theorem 3.5]). *Every nonnegative path and cycle flow $\hat{J}$ has a unique representation as a nonnegative arc flow $J$, given by* (13). *Furthermore, suppose $f : V \to \mathbb{R}$ satisfies $\sum_{v \in V} f_v = 0$, and suppose a nonnegative arc flow $J$ satisfies $D^\top J = f$. Then, $J$ can be decomposed as a path and cycle flow $\hat{J} : \mathcal{P} \cup \mathcal{C} \to \mathbb{R}_+$ such that every directed path with positive flow connects a source $s \in V$ with $f_s < 0$ to a target $t \in V$ with $f_t > 0$.*

We use this theorem to decompose the arc flow $J^\alpha$ into a path and cycle flow $\hat{J}^\alpha$. By optimality of $J^\alpha$, this decomposition does not contain any cycles $c \in \mathcal{C}$; if such a decomposition contained a cycle, then removing it can only decrease the total cost. We are left with path flows with positive mass.

Let $\mathcal{P}^{(s-t)}$ be the set of all directed paths connecting sources to targets:

$$\mathcal{P}^{(s-t)} = \left\{ r \in \mathcal{P} \mid f_{s(r)} < 0,\ f_{t(r)} > 0 \right\}.$$

By the remark above, for every $e \in E$,

$$J_e^\alpha = \sum_{r \in \mathcal{P}^{(s-t)}} \delta_e(r) \hat{J}^\alpha(r), \tag{14}$$

which we can denote in a more succinct notation $J^\alpha = \delta \cdot \hat{J}^\alpha$ where $\delta \in \{0,1\}^{|E| \times |\mathcal{P}^{(s-t)}|}$, is the matrix of path indicator functions. This decomposition of $J^\alpha$ into $(s-t)$ paths also applies to solutions of the non-regularized problem (LP).

For our remaining discussion, given a solution $J$ either of (LP) or (QP) for some $\alpha$, $\hat{J}$ will refer to a corresponding path flow above. Such a decomposition will be useful as it gives us a uniform upper bound on the flow depending only on $G$ and $f$:

**Proposition 5** (Boundedness). *Suppose $f : V \to \mathbb{R}$ satisfies $\sum_{v \in V} f(v) = 0$, and take a flow $J : E \to \mathbb{R}_+$ with $D^\top J = f$ decomposed into $\hat{J}$ as in (14). Then, for all paths $r \in \mathcal{P}^{(s-t)}$ and edges $e \in E$ we have $\hat{J}(r) \leq -f_{s(r)}$ and $J_e \leq -\sum_{v \in V, f_v < 0} f_v$.*

*Proof.* Since there are no cycles in the decomposition of $J$, there exists a source $s_1 \in V$ with no incoming flow. By the divergence constraint, $\hat{J}(r) \leq -f_{s_1}$ for all $r \in \mathcal{P}^{(s-t)}$ such that $s(r) = s_1$, i.e., the outgoing flow cannot be larger than what the source has to offer.

Remove this source $s_1$ and all the paths starting from $s_1$ from the graph, and repeat the process: There exists another source $s_2 \in V$ that has no incoming flow otherwise we would have a cycle, so by the same reasoning $\hat{J}(r) \leq -f_{s_2}$ for all $r \in \mathcal{P}^{(s-t)}$ such that $s(r) = s_2$. Repeating this process until exhausting all the sources shows $\hat{J}(r) \leq -f_{s(r)}$ for all $r \in \mathcal{P}^{(s-t)}$. Furthermore, for every edge $e \in E$:

$$J_e = \sum_{r \in \mathcal{P}^{(s-t)}} \delta_e(r)\hat{J}(r) = \sum_{\substack{v \in V \\ f_v < 0}} \sum_{\substack{r \in \mathcal{P}^{(s-t)} \\ s(r)=v}} \delta_e(r)\hat{J}(r) \leq \sum_{\substack{v \in V \\ f_v < 0}} \left( \sum_{\substack{r \in \mathcal{P}^{(s-t)} \\ s(r)=v}} \hat{J}(r) \right) = -\sum_{\substack{v \in V \\ f_v < 0}} f_v$$

□

We make a final remark before beginning the proofs of the propositions above. Given a nonnegative flow $J$ decomposed as $J = \delta \hat{J}$, we can rewrite the objective function $V_\alpha$ of (QP) on the path flows:

$$V_\alpha(J) = c^\top J + \frac{\alpha}{2}|J|^2 = \check{c}^\top \hat{J} + \frac{\alpha}{2}\hat{J}^\top S \hat{J} := \hat{V}_\alpha(\hat{J}), \tag{15}$$

where $\check{c} := \delta^\top c \in \mathbb{R}^{|\mathcal{P}^{(s-t)}| \times 1}$ and $S := \delta^\top \delta \in \mathbb{R}^{|\mathcal{P}^{(s-t)}| \times |\mathcal{P}^{(s-t)}|}$. Note that $S$ is a symmetric nonnegative matrix and $0 \leq S(r, r') \leq |E|$ for all $r, r' \in \mathcal{P}^{(s-t)}$.

### 4.4 Proofs of Proposition 4 and Corollary 1

Now we are ready to begin the proofs of proposition 4 and corollary 1. For the first proof, we will use the following lemma:

**Lemma 1** (Divergence-free decomposition). *Suppose $f : V \to \mathbb{R}$ satisfies $\sum_{v \in V} f_v = 0$, and take $J_1, J_2 : E \to \mathbb{R}_+$ with $D^\top J_1 = D^\top J_2 = f$; assume these flows admit decompositions $\hat{J}_1, \hat{J}_2$ with $J_i = \delta \hat{J}_i$. Then, there exist paths $X_-^1, X_+^1, \ldots, X_-^n, X_+^n \subseteq \mathcal{P}^{(s-t)}$ satisfying*

$$\begin{cases} \hat{J}_2(r) < \hat{J}_1(r) & \text{for all } r \in X_-^k, k \in \{1, \ldots, n\} \\ \hat{J}_2(r) > \hat{J}_1(r) & \text{for all } r \in X_+^k, k \in \{1, \ldots, n\} \end{cases}$$

*and corresponding scalars $\epsilon_1, \ldots, \epsilon_n > 0$ such that $\hat{J}_1 = \hat{J}_2 + \sum_{k=1}^n \epsilon_k \hat{R}_k$, where*

$$(R_k)_e := \sum_{r \in X_-^k} \delta_e(r) - \sum_{r \in X_+^k} \delta_e(r). \tag{16}$$

*The arc flows $R_k$ satisfy $D^\top R_k = 0$ for all $k \in \{0, \ldots, n\}$.*

Note we can equivalently write (16) via the path decomposition

$$\hat{R}_k = \sum_{r \in X_-^k} \chi_r - \sum_{r \in X_+^k} \chi_r,$$

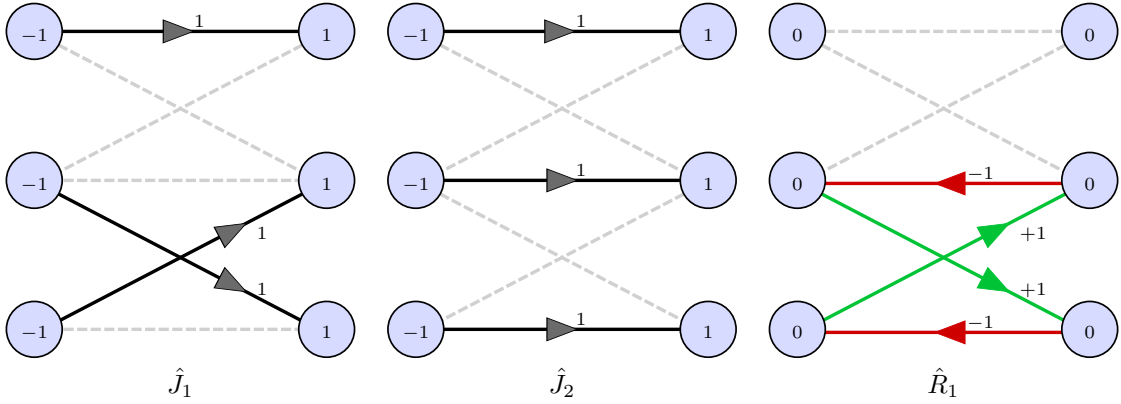

Figure 3: Example decomposition provided by Lemma 1. Here, $\hat{J}_1 = \hat{J}_2 + \varepsilon_1 \hat{R}_1$, where $\varepsilon = 1$.

where $\chi_r \in \{0,1\}^{|\mathcal{P}^{(s-t)}| \times 1}$ is the indicator of the $r$-th path. Here the $+/-$ subscript denotes excess and lack of mass, respectively. Figure 3 illustrates a simple example in which we need only one $\hat{R}_k$ term.

Our basic approach to proving this lemma will be to construct the $R_k$'s one-at-a-time so that $\hat{J}_1$ and $\hat{J}_2 + \sum_{k=1}^{\ell} \epsilon_\ell \hat{R}_\ell$ agree on at least one more path than in step $\ell-1$. If they do not agree completely, we identify a set of $(s-t)$ paths $X_+^k$ where $\hat{J}_1(r) < \hat{J}_2(r)$ and another set of $(s-t)$ paths $X_-^k$ where $\hat{J}_1(r) > \hat{J}_2(r)$. Furthermore, we enforce the following property for the sources and targets of both sets of $(s-t)$ paths: Each source has exactly one outgoing path from $X_-^k$ and one outgoing path from $X_+^k$, and each target has exactly one incoming path from $X_-^k$ and one incoming path from $X_+^k$. We repeat this procedure until all paths agree.

*Proof.* If $\hat{J}_1 = \hat{J}_2$ then choose $n = 0$, and there is nothing to prove.

Assume $\hat{J}_1 \neq \hat{J}_2$. Then, there exists a path $r_0 \in \mathcal{P}^{(s-t)}$ such that $\hat{J}_1(r_0) \neq \hat{J}_2(r_0)$; we can assume without loss of generality that $\hat{J}_1(r_0) > \hat{J}_2(r_0)$. This implies the source $s_0 = s(p_0)$ sends strictly more mass by $\hat{J}_1$ than $\hat{J}_2$ to the target $t_0 = t(r_0)$ along $r_0$. By the divergence constraint, another path $r_1 \in \mathcal{P}^{(s-t)}$ must have the opposite relationship, i.e. $\hat{J}_1(r_1) < \hat{J}_2(r_1)$. Given $r_0$, $s_0 = s(r_0)$ and $t_0 = t(r_0)$, we will iterate this basic observation to construct $\epsilon_k$ and paths $X_-^k, X_+^k$.

We phrase our proof in terms of the following algorithm. Assume we initialize $\hat{J}_2' := \hat{J}_2$, and loop over $m = 0, 1, 2, \ldots$ as follows:

- **Step 1:** Construct $X_-^{m+1}, X_+^{m+1}$:

  1. Since $\hat{J}_1 \neq \hat{J}_2'$, identify $r_0 \in \mathcal{P}^{(s-t)}$ such that $\hat{J}_1(r_0) \neq \hat{J}_2(r_0)$. Assume without loss of generality that $\hat{J}_1(r_0) > \hat{J}_2'(r_0)$.

  2. For $n = 1, 2, \ldots$, repeat the following, where $s_n := s(r_n)$ and $t_n := t(r_n)$:

     - If $n$ is odd, add a new target $t_n \in V$ and a path $r_n \in \mathcal{P}^{(s-t)}$ such that $s(r_n) = s_{n-1}$, $t(r_n) = t_n$, and $\hat{J}_2'(r_n) > \hat{J}_1(r_n)$.
       Such a path must exist: If it did not, then all paths $r$ leaving $s_{n-1}$ satisfy $\hat{J}_2'(r) \leq \hat{J}_1(r)$. Since $r_{n-1}$ leaving $s_{n-1}$ satisfies $\hat{J}_2'(r) < \hat{J}_1(r)$, this would contradict the divergence constraint:
       $$-f_{s_{n-1}} = \sum_{\substack{r \in \mathcal{P}^{(s-t)} \\ s(r) = s_{n-1}}} \hat{J}_2(r) < \sum_{\substack{r \in \mathcal{P}^{(s-t)} \\ s(r) = s_{n-1}}} \hat{J}_1(r) = -f_{s_{n-1}}.$$

     - Symmetrically, if $n$ is even, add a new source $s_n \in V$ and a path $r_n \in \mathcal{P}^{(s-t)}$ such that $s(r_n) = s_n$, $t(r_n) = t_{n-1}$, and $\hat{J}_2'(r_n) < \hat{J}_1(r_n)$.

- If $s_n$ or $t_n$ is a duplicate, halt this inner loop.

3. Depending on how we terminated the loop above, select a specific subgraph to construct $X_-^{m+1}$ and $X_+^{m+1}$:

    - Case 1, duplicated a source $s_l$ ($0 \leq l \leq n-2$, $l$ even), or a target $t_l$ ($0 < l \leq n-2$, $l$ odd) : Take

    $$X_-^{m+1} := \{r_j : l+1 \leq j \leq n, j \text{ even}\} \text{ and}$$
    $$X_+^{m+1} := \{r_j : l+1 \leq j \leq n, j \text{ odd}\}.$$

    - Case 2, repeated $t_0$: Take

    $$X_-^{m+1} := \{r_j : 0 \leq j \leq n, j \text{ even}\} \text{ and}$$
    $$X_+^{m+1} := \{r_j : 0 \leq j \leq n, j \text{ odd}\}.$$

    These two cases are constructed so that $X_+^{m+1}$ and $X_-^{m+1}$ satisfy the property: ($\star$) every source or target is connected with exactly one path in $X_-^{m+1}$ and one path in $X_+^{m+1}$, or is not connected at all.

- **Step 2:** Set

$$\hat{R}_{m+1} := \sum_{r \in X_-^{m+1}} \chi_r - \sum_{r \in X_+^{m+1}} \chi_r$$

The arc flow $R_{m+1}$ from $\hat{R}_{m+1}$ is divergence-free by ($\star$) above, i.e. $D^\top R = 0$. By construction, $\hat{R}_{m+1}$ removes mass where $\hat{J}_2 > \hat{J}_1$ and adds mass where $\hat{J}_1 > \hat{J}_2$.

- **Step 3:** Set

$$\epsilon_{m+1} := \min_{r \in X_-^{m+1} \cup X_+^{m+1}} |\hat{J}_1(r) - \hat{J}_2'(r)|$$

By construction, $\hat{J}_1$ and $\hat{J}_2$ never coincide on paths $r \in X_-^{m+1}$ and $X_+^{m+1}$, so $\epsilon > 0$.

- **Step 4:** Update $\hat{J}_2'$ via

$$\hat{J}_2' := \hat{J}_2 + \sum_{k=1}^{m+1} \epsilon_k \hat{R}_k.$$

$\hat{J}_2'$ is still nonnegative since we remove mass on the paths where $\hat{J}_1(r) < \hat{J}_2'(r)$, with an amount $\epsilon_{m+1}$ that maintains nonnegativity. Taking

$$r_{m+1}^* := \arg\min_{r \in X_-^{m+1} \cup X_+^{m+1}} |\hat{J}_1(r) - \hat{J}_2'(r)|,$$

by construction $\hat{J}_1(r_{m+1}^*) = \hat{J}_2'(r_{m+1}^*)$. The previously constructed $r_1^*, \ldots, r_m^* \notin X_-^{m+1} \cup X_+^{m+1}$, maintaining the property $\hat{J}_1(r_k^*) = \hat{J}_2'(r_k^*)$ for all $k \in \{1, \ldots, m\}$ after the update.

- **Step 5:** Increment $m$, and then halt if $\hat{J}_1 = \hat{J}_2'$. Otherwise, return to Step 1.

This algorithm must terminate because at each iteration we ensure that there is at least an extra path $r \in \mathcal{P}^{(s-t)}$ such that $\hat{J}_1(r) = \hat{J}_2'(r)$. At this point, $\hat{J}_1 = \hat{J}_2'$, or after expansion, $\hat{J}_1 = \hat{J}_2 + \sum_k \epsilon_k \hat{R}_k$. □ □

We use this lemma to prove Proposition 4:

*Proof of Proposition 4.* We prove the proposition by contradiction. To start, suppose $\alpha \in (0, \tilde{\alpha})$ for some $\tilde{\alpha} > 0$ that will be quantified later depending exclusively on the graph $G$ and the data $f$. Let $J^\alpha$ be the unique solution of (QP) and assume that $J^\alpha$ does not solve (LP).

Our proof leverages the fact that the $L^1$ cost dominates when $\alpha$ is small. So, if we move some mass onto some edges that are active in (LP), we can decrease the $L^1$ cost without significantly perturbing the $L^2$ term. Our argument proceeds in three steps:

1. Perturb the flow $J^\alpha$ into $J'$ by moving mass

2. Show that this perturbed flow $J'$ decreases the $L^1$ cost

3. Show that for $\alpha$ small enough, the perturbed flow performs better than $J^\alpha$ for (QP), which yields a contradiction.

The biggest difficulty of this proof is finding a relevant divergence-free perturbation of the flow; this task is facilitated by the decomposition in Lemma 1.

For a solution $J_0$ of (LP), Lemma 1 shows there exists $n \in \mathbb{N}$, $\epsilon_1, \ldots, \epsilon_n > 0$ and divergence-free arc flows $R_1, \ldots, R_n$ such that $\hat{J}^0 = \hat{J}^\alpha + \sum_{k=1}^n \epsilon_k \hat{R}_k$. Moreover, $n > 0$ since we assumed that $J^\alpha$ doesn't solve (LP). We begin by proving $\check{c}^\top \hat{R}_k \leq 0$ for all $k \in \{1, \ldots, n\}$; for ease of reading the vertical line on the left demarcates the proof of this claim.

Suppose there exists $k \in \{1, \ldots, n\}$ such that $\check{c}^\top \hat{R}_k > 0$. Define

$$\epsilon' := \frac{1}{2} \min_{r \in X_-^k} \{\hat{J}_0(r)\}.$$

Note $\epsilon' > 0$ since by construction in Lemma 1, $\hat{J}_0(r) > \hat{J}^\alpha(r) \geq 0$ for $r \in X_-^k$.

Next, define $\hat{J}'_0 := \hat{J}_0 - \epsilon' \hat{R}_k$. Consider a path $r$. Recalling $\hat{J}_0 \geq 0$, we have three cases:

$\boxed{r \in X_-^k}$ $\hat{J}_0(r) > \epsilon' \hat{R}_k(r)$, since $\hat{R}_k(r) = 1$ and by choice of $\varepsilon'$.

$\boxed{r \in X_+^k}$ $\hat{R}_k(r) = -1$, so again $\hat{J}'_0(r) \geq 0$ because $-\varepsilon' \hat{R}_k(r) > 0$.

$\boxed{r \notin X_-^k \cup X_+^k}$ $\hat{R}_k(r) = 0$ so again we still have $\hat{J}'_0(r) \geq 0$.

Hence, $\hat{J}'_0$ is nonnegative, and so is $J'_0$ after applying (13). Furthermore, since $D^\top R^k = 0$, we have $D^\top J'_0 = f$, implying that $J'_0$ is a feasible solution.

Finally,

$$c^\top J'_0 = \check{c}^\top \hat{J}'_0 = \check{c}^\top \hat{J}_0 - \epsilon' \check{c}^\top \hat{R}_k < \check{c}^\top \hat{J}_0 = c^\top J_0$$

since we assumed $\check{c}^\top \hat{R}_k > 0$. This contradicts optimality of $J_0$ for (LP), invalidating our initial assumption. Hence, $\check{c}^\top \hat{R}_k \leq 0$, $\forall k$. ◇

We next show that there exists some $k \in \{1, \ldots, n\}$ such that $\check{c}^\top \hat{R}_k < 0$:
Suppose $\check{c}^\top \hat{R}_k = 0$ for all $k \in \{1, \ldots, n\}$. Then,

$$c^\top J_0 = \check{c}^\top \hat{J}_0 = \check{c}^\top \hat{J}^\alpha + \sum_{k=1}^n \epsilon_k \check{c}^\top \hat{R}_k = \check{c}^\top \hat{J}^\alpha = c^\top J^\alpha$$

So, $J^\alpha$ solves (LP). This contradicts our assumption at the outset that $J^\alpha$ is *not* a solution of (LP).
◇

Now that we found which paths to perturb, we will use them to decrease the $L^1$ cost while judiciously controlling the increase of the $L^2$ term via an appropriate choice of constants. Introduce the perturbed flow $\hat{J}' = \hat{J}^\alpha + \epsilon \hat{R}_k$, where $\epsilon$ is a small constant that will be chosen later. This perturbed flow is nonnegative by construction of $X_-^k, X_+^k$ and satisfies $D^\top J' = f$. Our remaining task is to show that it performs better than $J^\alpha$ for (QP), or equivalently that the following value is smaller than $\hat{V}_\alpha(\hat{J})$:

$$\hat{V}_\alpha(\hat{J}') = \hat{V}_\alpha(\hat{J}^\alpha) + \epsilon \left[ \check{c}^\top \hat{R}_k + \alpha(\hat{J}^\alpha)^\top S \hat{R}_k \right] + \frac{\epsilon^2}{2} \alpha \hat{R}_k^\top S \hat{R}_k. \tag{17}$$

Define

$$K_{\min} := \max_{X_1, X_2 \subseteq \mathcal{P}^{(s-t)}} \left\{ \sum_{r \in X_1} \check{c}(r) - \sum_{r \in X_2} \check{c}(r) : \sum_{r \in X_1} \check{c}(r) - \sum_{r \in X_2} \check{c}(r) < 0 \right\} \tag{18}$$

$$= \max_{X_1, X_2 \subseteq \mathcal{P}^{(s-t)}} \left\{ \check{c}^\top \hat{R} : \hat{R} = \sum_{r \in X_1} \chi_r - \sum_{r \in X_2} \chi_r, \ \check{c}^\top \hat{R} < 0 \right\}.$$

Note that $K_{\min} < 0$ since $|\mathcal{P}^{(s-t)}| < +\infty$.

By Proposition 5, there exists a constant $M_1 > 0$ that only depends on the graph and $f$ such that $|\hat{J}^\alpha| \leq M_1$. Furthermore, since there is a finite number of paths in $G$, there exists $M_2 > 0$ which only depends on the graph such that $|\hat{R}_k| = |\sum_{r \in X_-} \chi_r - \sum_{r \in X_+} \chi_r| \leq M_2$. Finally, the matrix $S$ in (15) satisfies $|S| \leq M_3$ for some $M_3 > 0$ depending only on the graph, since every entry is bounded by $|E|$. $|S|$ denotes the spectral norm of $S$. Combining these inequalities, there exists $K_1 > 0$ that only depends on the graph and $f$ such that $|(\hat{J}^\alpha)^\top S \hat{R}_k| < K_1$. Choose $\tilde{\alpha} = \frac{|K_{\min}|}{2K_1} > 0$, and assume $0 < \alpha < \tilde{\alpha}$.

By construction, for such a choice of $\alpha$ we have

$$\begin{aligned}
\check{c}^\top \hat{R}_k + \alpha(\hat{J}^\alpha)^\top S \hat{R}_k &< \check{c}^\top \hat{R}_k + \alpha K_1 && \text{since } |(\hat{J}^\alpha)^\top S \hat{R}_k| < K_1 \\
&< \check{c}^\top \hat{R}_k + \frac{|K_{\min}|}{2} && \text{by choice of } \tilde{\alpha} \\
&< \frac{K_{\min}}{2} < 0 && \text{by definition of } K_{\min}.
\end{aligned} \tag{19}$$

Take $K_3 := \lambda_{\max}(S) M_2^2$, where $\lambda_{\max}(S)$ denotes the maximum eigenvalue of $S$; note $K_3$ only depends on the structure of the graph. This choice provides the inequality $\alpha \hat{R}_k^\top S \hat{R}_k < \tilde{\alpha} \lambda_{\max}(S) |\hat{R}_k|^2 \leq \tilde{\alpha} K_3$. Finally, choose $\epsilon > 0$ small enough such that $\hat{J}'$ remains nonnegative and $\frac{K_{\min}}{2} \epsilon + \tilde{\alpha} K_3 \frac{\epsilon^2}{2} < 0$; such an $\epsilon$ must exist because $K_{\min} < 0$. Then, by (17) and (19) we have $\hat{V}_\alpha(\hat{J}') < \hat{V}_\alpha(\hat{J}^\alpha)$. ⋄

This final inequality shows $V_\alpha(J') < V_\alpha(J^\alpha)$, which contradicts the definition of $J^\alpha$ as the unique minimizer of $V_\alpha$. □ □

The proof of Corollary 1 is rather straightforward:

*Proof of Corollary 1.* Pick $\tilde{\alpha}$ given by Proposition 4, and let $\alpha, \alpha'$ such that $0 < \alpha < \alpha' < \tilde{\alpha}$.

Define $J^\alpha$ (resp. $J^{\alpha'}$) to be the solution of the quadratically regularized problem ((QP)) with a regularization $\alpha$ (resp. $\alpha'$).

According to Proposition 4, $J^\alpha, J^{\alpha'}$ are also solutions of (LP). We will show, by contradiction, that they are the same solution.

Assume $J^\alpha \neq J^{\alpha'}$.

Because both flows are solutions of (LP), we have that $\check{c}^\top \hat{J}^\alpha = \check{c}^\top \hat{J}^{\alpha'}$.

Since $J^\alpha$ is the unique solution of (QP),

$$V_\alpha(J^\alpha) = \hat{V}_\alpha(\hat{J}^\alpha) = \check{c}^\top \hat{J}^\alpha + \frac{\alpha}{2}(\hat{J}^\alpha)^\top S \hat{J}^\alpha < \check{c}^\top \hat{J}^{\alpha'} + \frac{\alpha}{2}(\hat{J}^{\alpha'})^\top S \hat{J}^{\alpha'} = V_\alpha(J^{\alpha'}).$$

because $\check{c}^\top \hat{J}^\alpha = \check{c}^\top \hat{J}^{\alpha'}$, we have that $(\hat{J}^\alpha)^\top S \hat{J}^\alpha < (\hat{J}^{\alpha'})^\top S \hat{J}^{\alpha'}$, contradicting Proposition 3.

We conclude that $J^\alpha = J^{\alpha'}$ is the same solution of LP. □ □

## 4.5 Monotonicity conjecture

As shown in a previous counterexample, we do not necessarily have monotonicity of the active sets $S(\alpha)$ as $\alpha$ grows. Indeed, when $\alpha$ increases, new "$L^2$-friendly" edges might be activated thanks to a mass transfer from previously active "$L^1$-friendly" edges. In certain cases, this mass transfer completely deactivates these $L^1$-friendly edges, ruling out such a property.

We experimentally observe a monotonicity property in the mass flow on the edges, however. It appears that mass flow on these $L^2$-friendly edges tends to always grow as $\alpha$ increases, and vice versa for the $L^1$-friendly edges. One difficulty consists in determining these edges, which depends in a non-trivial way on the graph $G$ and the cost $c$.

This monotonicity behavior is summarized in the following conjecture:

**Conjecture 1** (Monotonicity). *Let $J^{\alpha_1}, J^{\alpha_2}$ be the solutions of the quadratic program* (QP) *for $0 < \alpha_1 < \alpha_2$. Use Lemma 1 to write*

$$J^{\alpha_2} = J^{\alpha_1} + \sum_i \epsilon_i R_i,$$

*where $\epsilon_i > 0$. Then, for all $\alpha' \in (\alpha_1, \alpha_2)$, the solution $J^{\alpha'}$ of the quadratic program* (QP) *with regularization $\alpha'$ is given by*

$$J^{\alpha'} = J^{\alpha_1} + \sum_i \epsilon'_i R_i$$

*with the same divergence free flows $R_i$, for some $\epsilon'_i \in [0, \epsilon_i]$.*

If this is true, one can choose the parameter $\alpha_2 \to +\infty$, so that the minimization problem (QP) would yield a minimizer $J^\infty$ for the quadratic term only, regardless of the cost $c$. Comparing this solution to $\lim_{\alpha_1 \to 0} J^{\alpha_1} = J^0$, and in particular extracting the divergence-free flows $R_i$ from the difference $J^\infty - J^0$, yields a complete characterization of the $L^1$- and $L^2$-friendly paths. This result still does not give a complete description of $\epsilon_i$'s in terms of $\alpha$, but we believe that as it is stated, Lemma 1 should be enough to obtain a clean proof.

Figure 4 illustrates the monotonicity behavior on a graph with 5 vertices and 7 edges. The cost per edge is 1 unit per mass flow for all edges, except for the edge leaving the source and pointing obliquely to the right which has cost 10. The red vertex is a source providing 1 unit of mass (negative mass), the green vertex is a sink taking away 1 unit of mass (positive); all other vertices are intermediary nodes.

In the top three parts of Figure 4, values of the flow solving (QP) for $\alpha = 10^{-1}, 10$ and $10^3$ are marked on the edges, with positive flow represented in green, and negative flow in red. In the bottom two parts, values of the divergence-free flows $R_i$, where $J^\alpha = J^0 + \sum_i \epsilon_i R_i$, are represented.

When $\alpha = 10^{-1}$, we obtain a solution of LP. As $\alpha$ increases to 10, the divergence-free loop $R_1$ illustrated below the case $\alpha = 10$ in Figure 4 is activated, with a coefficient growing up to $\epsilon_1 = 0.07$. Increasing $\alpha$ will further increase $\epsilon_1$, meaning that mass flow on the newly activated edges will increase and mass flow on a previously activated edge will decrease. A similar behavior occurs when $\alpha$ varies from 10 to $10^3$; on top of $\epsilon_1 = 0.25$, another divergence-free loop $R_2$ gets activated with a coefficient $\epsilon_2 = 0.25$.

## 5 Optimization Algorithm

Our optimization algorithm for quadratically-regularized flow takes inspiration from classical iterative methods for solving linear systems. In each iteration, we choose a search direction $s$ from the current iterate $x$ and update our iteration via line search $x \mapsto x + ts$ for some $t > 0$; note that this standard approach from optimization actually contrasts with the alternating projection algorithms commonly used e.g. in [3] for regularized transport. The difference here is that our piecewise-quadratic objective function makes it possible to carry out line search in closed-form.

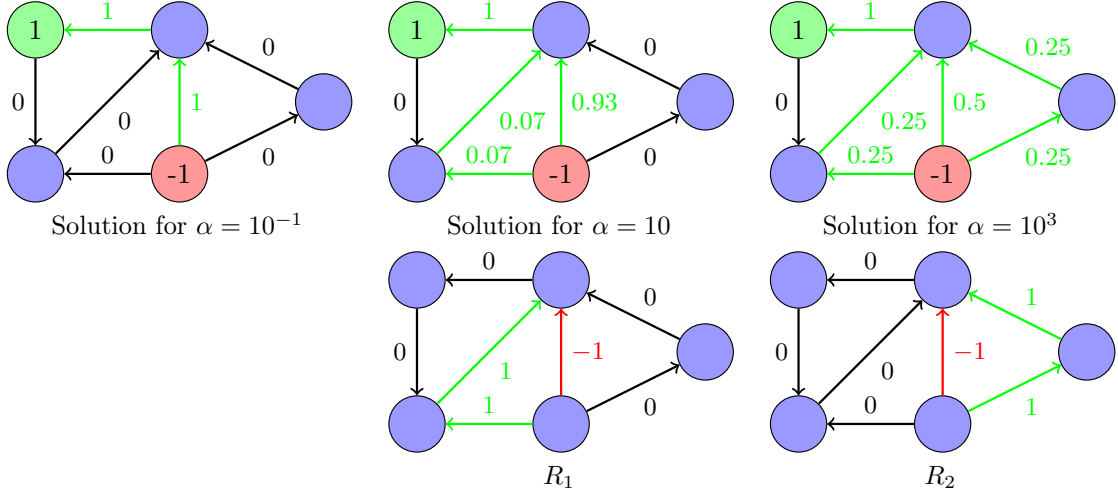

Figure 4: Example of solutions of (QP) for increasing values of $\alpha$, and the corresponding newly activated divergence free flows $R_i$

From a very high level, our iterative algorithm for optimizing the dual problem (5) divides into standard steps:

$$s_k \leftarrow \text{SEARCH-DIRECTION}(p_k) \qquad \S 5.2$$
$$t_k \leftarrow \text{LINE-SEARCH}(p_k, s_k) \qquad \S 5.1$$
$$p_{k+1} \leftarrow p_k + t_k s_k. \tag{20}$$

We will present these steps in reverse order, since our line search procedure will inform our choice of search directions.

## 5.1 Line Search

Our first goal is to define a procedure for increasing the dual objective (5) given a current iterate $p_k \in \mathbb{R}^{|V|}$ and a search direction $s_k \in \mathbb{R}^{|V|}$. Denote $M_k := M(p_k) \in \{0,1\}^{|E| \times |E|}$, the diagonal matrix indicating the active set $S(p_k)$. Then, for small $t > 0$, $p_k$ the dual objective restricted to the line through $s_k$ is the parabola

$$g_k(t) := \alpha f^\top (p_k + ts_k) - \frac{1}{2}(Dp_k + tDs_k - c)^\top M_k (Dp_k + tDs_k - c)$$
$$= -\frac{1}{2}t^2[s_k^\top L_k s_k] + t[\alpha f^\top s_k - v_k M_k D s_k] + \text{const.},$$

where $L_k := D^\top M_k D$ is the Laplacian of the active subgraph and $v_k := Dp_k - c$. This parabola is minimized at

$$t_{\text{quadratic}} := \frac{\alpha f^\top s_k - v_k^\top M_k D s_k}{s_k^\top L_k s_k}. \tag{21}$$

Note $t_{\text{quadratic}} \geq 0$ in every step of our algorithm because we choose $s_k$ to be an ascent direction for the dual objective.

This formula for minimizing the current parabola, however, is only applicable if $M(p_k) = M(p_k + t_{\text{quadratic}} s_k)$, that is, while the active set of edges $S(\cdot)$ remains unchanged. Hence, we limit our line search if $M$ changes before the parabola is minimized. Define the "hitting time" vector

$$h_k := -v_k \oslash (Ds_k), \tag{22}$$

where $\oslash$ indicates elementwise division. Then, if we define

$$t_{\text{active set}} := \begin{cases} \min_{e \in E} & h_{ke} \\ \quad \text{s.t.} & h_{ke} > 0, \end{cases} \quad (23)$$

then we choose

$$t_k \leftarrow \min(t_{\text{quadratic}}, \ t_{\text{active set}}). \quad (24)$$

If $t_k = t_{\text{active set}}$, then the active set and hence $M_k$ changes. In §5.2.2, we document how to deal with this change elegantly. We never find a case where $t_{\text{active set}} = 0$ before convergence.

## 5.2 Search Direction

All that remains is the choice of search direction $s_k$ given the current iterate $p_k$. Ideally, we might wish to use a search direction from Newton's method, but the Hessian (9) potentially has a large null space of dimensionality equal to the number of connected components in the active graph.

Locally, the objective is quadratic but not strictly convex; in particular, there are directions along which it may be flat. Hence, we adopt an *alternating* strategy:

$$s_k \leftarrow \begin{cases} \alpha f - D^\top M_k v_k & \text{if } k \text{ is odd} \\ L_k^+(\alpha f - D^\top M_k v_k) & \text{if } k \text{ is even} \end{cases} \quad \begin{array}{l} \textit{Gradient direction} \\ \textit{Pseudo-Newton direction.} \end{array} \quad (25)$$

The objective (5) is not affected by adding any multiple of $\mathbb{1}$ to $p$, so for numerical stability we shift $s_k$ to sum to zero. We use the phrase "pseudo-Newton" to refer to the fact that the Hessian $L_k$ is not invertible. The Moore–Penrose pseudoinverse $L_k^+$ essentially restricts optimization to the space in which the objective function has curvature; then, the gradient step accounts for optimization in directions for which the Hessian provides no information.

The experiments in §6 verify that our choice of search directions is useful in practice for reducing the *number* of iterations. But, as written it appears that every iteration of our algorithm is far more expensive algorithmically: The Laplacian $L_k$ changes any time the active set changes, preventing the use of a fixed factorization to apply $L_k^+$. Our experiments also show that preconditioning gradient descent using the Laplacian $L := D^\top D$ of the entire graph $G$ is ineffective. Hence, as-written the algorithm could take upwards of $O(|V|^3)$ time to apply $L_k^+$, whereas individual iterations of gradient descent take $O(|V|)$ time.

As suggested in [29], however, adding or removing an edge to the active set $S(p_k)$ corresponds to a **rank-1 change** of the Laplacian $L_k$. In particular, recall that $L_k := D^\top M_k D$, where $M_k$ contains 1 on the diagonal whenever the corresponding edge is in the active set. Adding an edge to the active set corresponds to flipping a diagonal element of $M_k$ from 0 to 1, while removing an edge does the opposite. In other words, if $t_k = t_{\text{active set}}$, then $L_{k+1} = L_k \pm d_k d_k^\top$, where $d_k$ is the row of $D$ corresponding to the edge activated/deactivated at $t_{\text{active set}}$; otherwise $L_{k+1} = L_k$. This suggests that pre-factorization plus rank-1 updates similar to the Sherman–Morrison formula [37] can be used to make applying $L_k$ less expensive from iteration to iteration.

Two confounding factors, however, require specialized attention before we can use rank-1 update formulae:

1. Graph Laplacians are not full-rank, and hence standard techniques for Cholesky factorization will fail.

2. The pseudoinverse $L_k^+$ is a dense matrix, whereas $L_k$ is potentially sparse when $|E| \ll |V|^2$. Hence, rank-1 update formulas for $L_k^+$ directly as explained in [29] are still expensive computationally.

We develop strategies for overcoming these issues in the subsections below.

### 5.2.1 Dealing with Null Space

Tools for sparse matrix factorization and rank-1 updates typically assume that the matrices involved have full rank. In particular, they typically are optimized for applying inverses rather than pseudoinverses of matrices. Hence, we first provide a technique for writing application of $L_k^+$ in a sparse but low-rank fashion.

To do so, we note that the null space of $L_k$ is simply the set of indicators of each connected component in the active graph. Define $N_k \in \mathbb{R}_+^{|V| \times (\# \text{ components})}$ to be the matrix whose columns are an orthogonal basis for the null space of $L_k$, constructed e.g. by flood fill on the connected components of the active graph followed by normalization. Furthermore, define $P_k$ to be the projection operator onto the orthogonal complement of the column space of $N_k$:

$$P_k := I - N_k N_k^\top. \tag{26}$$

Suppose we wish to find the product $x := L_k^+ b$ for some $b \in \mathbb{R}^{|V|}$. Then, by definition of the pseudoinverse, $x$ satisfies the relationship $L_k x = Pb$. Furthermore, the pseudoinverse zeros out components of $x$ in the null space of $L_k$, showing $N_k^\top x = 0$, which of course implies $N_k N_k^\top x = 0$. Adding these two expressions together shows $(L_k + N_k N_k^\top)x = P_k b \implies x = (L_k + N_k N_k^\top)^{-1} P_k b$. In other words:

$$L_k^+ = (L_k + N_k N_k^\top)^{-1} P_k. \tag{27}$$

An analog of this formula also appears in [29].

### 5.2.2 Maintaining Sparsity

The matrix $L_k + N_k N_k^\top$ is invertible because the column spaces of $L_k$ and $N_k$ together span $\mathbb{R}^{|V|}$ by construction and hence can be factored using the standard Cholesky method. Unlike $L_k$, however, it is almost certainly dense.

To circumvent this issue, note:

$$L_k + N_k N_k^\top = \begin{pmatrix} M_k D \\ N_k^\top \end{pmatrix}^\top \underbrace{\begin{pmatrix} M_k D \\ N_k^\top \end{pmatrix}}_{W_k}. \tag{28}$$

Here, we leverage the idempotence property $M_k^2 = M_k$, since the diagonal of $M_k$ is composed of zeros and ones.

We apply sparse matrix QR factorization machinery proposed in [12] to factor $W_0 = Q_0 R_0$, providing an initial Cholesky factorization

$$L_0 + N_0 N_0^\top = W_0^\top W_0 = R_0^\top R_0.$$

Updating $R_k$ to $R_{k+1}$ is carried out by sparse rank-1 updates to the Cholesky factorization, illustrated in Figure 5:

- Introducing an edge $e \in E$ to the active set makes at most four rank-1 changes to $L + NN^\top$
  1. Adding $d_e d_e^\top$, where $d_e$ is the row of $D$ corresponding to $e$
  2. Subtracting $n_1 n_1^\top + n_2 n_2^\top$ and adding $\bar{n}\bar{n}^\top$, if $e$ merges two connected components $n_1$ and $n_2$ into a new component $\bar{n}$

- Removing an edge $e \in E$ from the active set makes at most four rank-1 changes to $L + NN^\top$
  1. Subtracting $d_e d_e^\top$, where $d_e$ is the row of $D$ corresponding to $e$
  2. Adding $n_1 n_1^\top + n_2 n_2^\top$ and subtracting $\bar{n}\bar{n}^\top$, if breaks a connected component $\bar{n}$ into two connected components $n_1$ and $n_2$

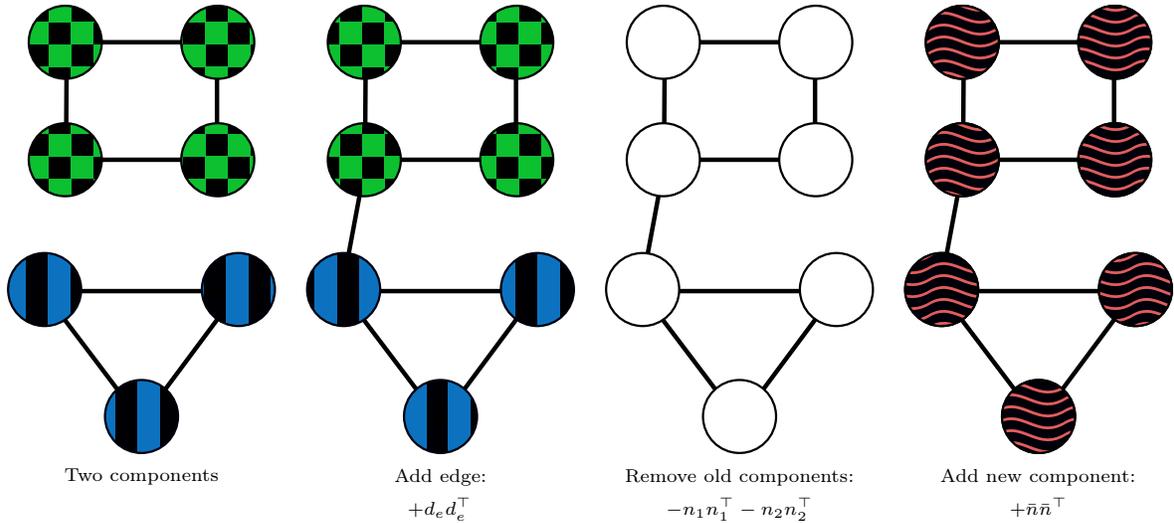

Figure 5: Rank-1 updates required to add an edge connecting two components. Removing an edge essentially corresponds to these operations in reverse.

Our algorithm never computes $L_k + N_k N_k^\top$ explicitly but rather updates its factorization $L_k + N_k N_k^\top = R_k R_k^\top$ using a sparse upper-triangular $R_k$ in each iteration via the rules above. Each rank-1 update is carried out using the CHOLMOD library in SuiteSparse, which uses the algorithm described in [7] to update/downdate a Cholesky factorization while maintaining sparsity.

With this rank-1 update machinery in place, Algorithm 1 provides pseudocode combining the steps we have described into a full algorithm for regularized graph transport.

## 6 Experiments

In this section, we illustrate the performance of our algorithm on graphs of different sizes, as well as for different regularization parameters $\alpha$. Our implementation is in Matlab, and rank-1 updates are dealt with using the library SuiteSparse. Our single-threaded implementation is run on a laptop with an Intel i7-4510U with 2 physical cores (4 logical) at 2.00GHz and 12 GB of RAM.

Our code, which we will refer to as HessUpdate, is benchmarked against two other codes; the first one, referred as GradDescent, is a classical gradient descent (or rather "ascent" in our case). The second one, referred as PrecondGrad, is similar to HessUpdate but does not update the Laplacian at all; it uses the full Laplacian of the graph at each Newton step, which is pre-factored initially.

The three implementations were benchmarked on 10 different graphs for each size, with 5 different sizes: 50, 100, 500, 1,000 and 5,000 nodes. These graphs were generated using the algorithm detailed in [41]. In particular, we used their method to generate random connected graphs of varying sizes whose vertex degrees are restricted to the range $[1, 10]$ and follow a heavy-tailed distribution with exponent $a = 2.5$ and average $z = 5$. For each graph, eight different parameters $\alpha$ are tested, $\alpha = 10^{-k}$ for $k \in \{-1, \ldots, 5\}$ as well as $\alpha = 5$ and $\alpha = 10$. For each graph and parameter $\alpha$, 10 runs are made with randomly generated data $f$. To generate the data $f$, we randomly select 10% of its nodes to be sources or sinks, and randomly assign values for their capacities uniformly chosen between $-10$ and $10$. Then we modify the capacity of the last node to be the opposite of the sum of all the other capacities, ensuring that $\sum_{v \in V} f_v = 0$.

The graphs generated are bidirectional with one connected component, which ensures that a feasible solution always exists for such random data. Results of these 10 runs are then averaged for each pair $(n_{\text{vertices}}, \alpha)$. The three methods are limited to a maximum of 3000 iterations, with converged

```
function REGULARIZED-DUAL-TRANSPORT(G = (V, E), c, α, f)
    // Optimization algorithm for dual problem (5)
    p ← RANDOM(|V| × 1)                                              // Randomly initialize dual variable p
    (v, M) ← Dp − c;   M ← diag(v > 0)                               // Indicator of active edges
    L ← D⊤MD                                                         // Laplacian of active subgraph
    N ← FLOOD-FILL-BASIS(L)                                          // Null space basis; one column per connected component
    R ← SPARSE-QR(((MD)⊤  N)⊤)                                       // Q not needed from QR factorization from (28)
    for k = 1, 2, 3, . . .
        // Choose line search direction
        s ← αf − D⊤Mv                                                // Gradient direction from (8)
        if k is even then                                            // Use pseudo-Newton every other iteration
            s ← R⁻¹R⁻⊤(s − NN⊤s)                                     // Forward then back substitution to apply L⁺ via (27)
        // Line search
        t_quadratic ← (αf⊤s − v⊤MDs) / (s⊤Ls)                        // Minimum of parabola, (21)
        h ← −v ⊘ (Ds)                                                // Time t at which each edge flips active or inactive, (22)
        t_active set ← min(POSITIVE-ELEMENTS(h))                     // First edge to change sign, (23)
        t ← min(t_quadratic, t_active set)                           // Go to minimum or to active set change (24)
        p ← p + ts                                                   // Updated dual variable, (20)
        // Update Cholesky factorization R⊤R = L + NN⊤ for new Laplacian
        (v, M) ← Dp − c;   M ← diag(v > 0)                           // Update indicator of active edges
        for each new edge e in active set
            L ← L + d_e d_e⊤                                         // Rank-1 addition to Laplacian matrix
            RANK-1-CHOLESKY-UPDATE(R, d_e)                           // Add d_e d_e⊤ to factorized matrix
            if adding e merges components c_1, c_2 into c̄ then
                REMOVE-COLUMNS(N, (c_1, c_2))                        // Remove old disjoint connected components
                ADD-COLUMN(N, c̄)                                     // Add new component
                RANK-1-CHOLESKY-UPDATE(R, c̄)                         // Add c̄c̄⊤ to factorized matrix
                RANK-1-CHOLESKY-DOWNDATE(R, c_1)                     // Subtract c_1 c_1⊤ from factorized matrix
                RANK-1-CHOLESKY-DOWNDATE(R, c_2)                     // Subtract c_2 c_2⊤ from factorized matrix
        for each edge e removed from active set
            L ← L − d_e d_e⊤                                         // Rank-1 subtraction from Laplacian matrix
            if removing e splits component c̄ into components c_1, c_2 then
                REMOVE-COLUMN(N, c̄)                                  // Remove old connected component
                ADD-COLUMNS(N, (c_1, c_2))                           // Add new individual components
                RANK-1-CHOLESKY-UPDATE(R, c_1)                       // Add c_1 c_1⊤ to factorized matrix
                RANK-1-CHOLESKY-UPDATE(R, c_2)                       // Add c_2 c_2⊤ to factorized matrix
                RANK-1-CHOLESKY-DOWNDATE(R, c̄)                       // Subtract c̄c̄⊤ from factorized matrix
            RANK-1-CHOLESKY-DOWNDATE(R, d_e)                         // Subtract d_e d_e⊤ from factorized matrix
    return p
```

Algorithm 1: Optimization algorithm for quadratically-regularized transport on graphs. See the accompanying MATLAB implementation for code including tolerances and convergence criteria. Cholesky factorization updates (adding $xx^\top$ for a vector $x$) are always scheduled before downdates (subtracting $xx^\top$) so that $R$ is always full-rank.

considered as achieved if the norm of the gradient is less than $10^{-8}$.

As a reference, we also include the runtime of a classical simplex algorithm on the unregularized problem, for the same graphs and initial datum. This algorithm will be referred as SIMPLEX, and uses MATLAB's dual simplex implementation via the command LINPROG.

Figure 6 provides the average run times of the code for each pair $(n_{\text{vertices}}, \alpha)$. Since this time is either the time to convergence or the time to reach the maximum number of allowed iterations, we will also compute the average percent error between the objective value when the code stops and the true global maximum value. These results are presented in Figure 7. Run times of SIMPLEX were independent of $\alpha$ (although each run was made with a different $f$); hence we only represent these results once, as a function of the size of the graph.

We also provide in the last table of Figure 7 the relative difference of the $L^1$ objective value $(c^T J)$ between the regularized solution with parameter $\alpha$ obtained by HESSUPDATE run until convergence and the unregularized solution obtained by SIMPLEX, for each run.

We can see from Figure 6 that HESSUPDATE always outperforms GRADDESCENT time-wise. This is unsurprising as checking every other iteration if a Newton direction is available allows us to bypass multiple gradient descent steps. We can also see from Figure 7 that convergence always occurs for HessUpdate within the allotted number of iterations except for high temperatures $\alpha \in \{5, 10\}$ on big graphs of size $1,000$ and $5,000$, whereas GradDescent always reaches the maximum number of iterations before converging, except for $\alpha \in \{10^{-3}, 10^{-2}\}$ on small graphs of size 50 and 100.

| Size | $\alpha = 10^{-5}$ | $10^{-4}$ | $10^{-3}$ | $10^{-2}$ | $10^{-1}$ | 1 | 5 | 10 |
|------|------|------|------|------|------|------|------|------|
| HessUpdate | | | | | | | | |
| 50   | $4.4 \cdot 10^{-3}$ | $6.48 \cdot 10^{-3}$ | $4.93 \cdot 10^{-3}$ | $5.73 \cdot 10^{-3}$ | $1.01 \cdot 10^{-2}$ | $1.73 \cdot 10^{-2}$ | $2.33 \cdot 10^{-2}$ | $3.42 \cdot 10^{-2}$ |
| 100  | $1.52 \cdot 10^{-2}$ | $1.16 \cdot 10^{-2}$ | $1.2 \cdot 10^{-2}$ | $1.39 \cdot 10^{-2}$ | $6.52 \cdot 10^{-2}$ | $7.51 \cdot 10^{-2}$ | 0.14 | 0.14 |
| 500  | 0.33 | 0.39 | 0.87 | 2.97 | 6.5 | 2.75 | 3.59 | 4.25 |
| 1000 | 1.2 | 6.32 | 6.78 | 8.6 | 11.14 | 13.61 | 13.77 | 13.55 |
| 5000 | 3.4 | 4.66 | 7.84 | 9 | 11.56 | 14.28 | 13.3 | 14.73 |
| GradDescent | | | | | | | | |
| 50   | 0.26 | 0.26 | 0.26 | 0.26 | 0.26 | 0.26 | 0.26 | 0.26 |
| 100  | 0.89 | 0.89 | 0.9 | 0.9 | 0.9 | 0.89 | 0.9 | 0.89 |
| 500  | 21.19 | 21.17 | 21.23 | 21.14 | 21.16 | 21.16 | 21.17 | 21.14 |
| 1000 | 87.88 | 87.86 | 87.93 | 87.89 | 87.9 | 88.05 | 88.02 | 87.88 |
| 5000 | 586.73 | 583.2 | 582.73 | 540.24 | 539.91 | 540.64 | 539.58 | 540.17 |
| PrecondGrad | | | | | | | | |
| 50   | 0.34 | 0.33 | 0.23 | $6.92 \cdot 10^{-2}$ | $2.18 \cdot 10^{-2}$ | $1.09 \cdot 10^{-2}$ | $6.41 \cdot 10^{-3}$ | $4.74 \cdot 10^{-3}$ |
| 100  | 1.77 | 1.7 | 1.52 | 0.62 | 0.25 | $8.57 \cdot 10^{-2}$ | $4.95 \cdot 10^{-2}$ | $3.67 \cdot 10^{-2}$ |
| 500  | 23.81 | 23.73 | 23.89 | 22.28 | 6.86 | 2.28 | 1.21 | 0.5 |
| 1000 | 94.67 | 94.65 | 94.87 | 94.53 | 45.2 | 29.56 | 8.5 | 4.04 |
| 5000 | 646.16 | 647.34 | 645.69 | 599.40 | 554.84 | 383.79 | 241.55 | 64.19 |
| Simplex *on the unregularized problem* | | | | | | | | |
| 50   | $1.07 \cdot 10^{-2}$ | | | | | | | |
| 100  | $1.44 \cdot 10^{-2}$ | | | | | | | |
| 500  | 0.21 | | | | | | | |
| 1000 | 2.66 | | | | | | | |
| 5000 | 14.38 | | | | | | | |

Figure 6: Average run times in seconds for HessUpdate, GradDescent and PrecondGrad function of the size of the graph and the parameter $\alpha$, and average run times for Simplex function of the graph size.

| Size | $\alpha = 10^{-5}$ | $10^{-4}$ | $10^{-3}$ | $10^{-2}$ | $10^{-1}$ | 1 | 5 | 10 |
|---|---|---|---|---|---|---|---|---|
| HessUpdate | | | | | | | | |
| 50   | 0               | 0               | 0               | 0               | 0               | 0               | 0               | 0 |
| 100  | $1.2 \cdot 10^{-3}$ | 0               | 0               | 0               | 0               | 0               | $4.8 \cdot 10^{-3}$ | $4 \cdot 10^{-4}$ |
| 500  | $1.4 \cdot 10^{-3}$ | 0               | $5 \cdot 10^{-4}$ | $1.87 \cdot 10^{-2}$ | $2.19 \cdot 10^{-2}$ | 0 | $3.84 \cdot 10^{-2}$ | $6.64 \cdot 10^{-2}$ |
| 1000 | $3.5 \cdot 10^{-3}$ | $1.5 \cdot 10^{-3}$ | $6.3 \cdot 10^{-3}$ | $3.13 \cdot 10^{-2}$ | $8.87 \cdot 10^{-2}$ | $2.93 \cdot 10^{-2}$ | 0.27 | 0.46 |
| 5000 | $5.9 \cdot 10^{-3}$ | $5.5 \cdot 10^{-3}$ | $1.23 \cdot 10^{-2}$ | $7.55 \cdot 10^{-2}$ | $8.9 \cdot 10^{-2}$ | 0.13 | 0.19 | 0.22 |
| GradDescent | | | | | | | | |
| 50   | 0.88 | 0.26 | $1.51 \cdot 10^{-2}$ | $1.55 \cdot 10^{-2}$ | $\infty$ | $\infty$ | $\infty$ | $\infty$ |
| 100  | 0.86 | 0.27 | $1.21 \cdot 10^{-2}$ | $2.02 \cdot 10^{-2}$ | $\infty$ | $\infty$ | $\infty$ | $\infty$ |
| 500  | 0.85 | 0.3 | $1.7 \cdot 10^{-2}$ | $\infty$ | $\infty$ | $\infty$ | $\infty$ | $\infty$ |
| 1000 | 0.87 | 0.32 | $1.82 \cdot 10^{-2}$ | $\infty$ | $\infty$ | $\infty$ | $\infty$ | $\infty$ |
| 5000 | 0.25 | $9.71 \cdot 10^{-2}$ | $4.5 \cdot 10^{-3}$ | $\infty$ | $\infty$ | $\infty$ | $\infty$ | $\infty$ |
| PrecondGrad | | | | | | | | |
| 50   | 0.22 | $7.09 \cdot 10^{-2}$ | $2.5 \cdot 10^{-3}$ | 0 | 0 | 0 | 0 | 0 |
| 100  | 0.26 | $7.48 \cdot 10^{-2}$ | $2.8 \cdot 10^{-3}$ | 0 | 0 | 0 | 0 | 0 |
| 500  | 0.37 | 0.1020 | $5.7 \cdot 10^{-3}$ | $1 \cdot 10^{-4}$ | 0 | 0 | 0 | 0 |
| 1000 | 0.43 | 0.1112 | $6.9 \cdot 10^{-3}$ | $1 \cdot 10^{-4}$ | 0 | 0 | 0 | 0 |
| 5000 | 0.52 | 0.3412 | $2.01 \cdot 10^{-2}$ | $1.2 \cdot 10^{-3}$ | $1 \cdot 10^{-4}$ | $1 \cdot 10^{-4}$ | 0 | 0 |
| *Relative difference in $L^1$ cost between* Simplex *and* HessUpdate | | | | | | | | |
| 50   | 0 | 0 | 0 | 0 | $1.35 \cdot 10^{-2}$ | $9.04 \cdot 10^{-2}$ | $6.12 \cdot 10^{-2}$ | $3.99 \cdot 10^{-2}$ |
| 100  | 0 | 0 | 0 | 0 | $1.4 \cdot 10^{-2}$ | $8.73 \cdot 10^{-2}$ | $5.56 \cdot 10^{-2}$ | $3.53 \cdot 10^{-2}$ |
| 500  | 0 | 0 | 0 | 0 | $1.9 \cdot 10^{-2}$ | $8.03 \cdot 10^{-2}$ | $4.84 \cdot 10^{-2}$ | $3.03 \cdot 10^{-2}$ |
| 1000 | 0 | 0 | 0 | 0 | $2 \cdot 10^{-2}$ | $7.41 \cdot 10^{-2}$ | $4.31 \cdot 10^{-2}$ | $2.67 \cdot 10^{-2}$ |
| 5000 | 0 | 0 | 0 | 0 | $2.15 \cdot 10^{-2}$ | $7.33 \cdot 10^{-2}$ | $4.23 \cdot 10^{-2}$ | $2.62 \cdot 10^{-2}$ |

Figure 7: Average relative error for HessUpdate, GradDescent and PrecondGrad as a function of the size of the graph and the parameter $\alpha$. $\infty$ means that the error was too big and the algorithm did not converge. Also represented is the average relative difference $L^1$ cost between HessUpdate and Simplex.

We note some interesting behavior when comparing HessUpdate and PrecondGrad; HessUpdate always outperforms in speed PrecondGrad for smaller regularizers $\alpha$ whereas the opposite is true for larger temperatures, as we can see from Figures 6 and 7. The explanation for this behavior is quite simple; our graphs have one connected component and are bidirectional. As we increase the temperature, the active graph tends to grow, barring pathological counterexamples. We can experimentally check that for high temperatures, the optimal flow will activate nearly all edges on the tested graphs, and hence the full Laplacian is a good approximation for the active Laplacian of the optimal solution. This allows PrecondGrad to have quicker iteration because there is no update to do for the Laplacian, while conserving good convergence speed because of a good static guess for the active graph.

For smaller temperatures, the opposite is true; the full Laplacian is a bad approximation as we can guess from the sparsity property 4. Hence the advantage gained from not having to compute or update the Laplacian at each step is minimal compared to the loss of speed due to slow convergence.

The comparison between HessUpdate and Simplex from Figure 6 shows that overall computations times of the regularized and unregularized problem are in the same ballpark, although the simplex algorithm has some scaling issues. Figure 7 shows that HessUpdate's solutions for the regularized problem are actually the same as for the unregularized problem for temperatures $\alpha \leq 10^{-2}$. This is consistent with the sparsity proposition 4. As $\alpha$ grows, however, the regularized solution starts greatly differing from the unregularized one.

| Grid Size | $\alpha = 10^{-4}$ | $10^{-3}$ | $10^{-2}$ | $10^{-1}$ | 1 | 5 | 10 |
|---|---|---|---|---|---|---|---|
| HessUpdate | | | | | | | |
| 10 | 0.9 | 0.11 | $5.84 \cdot 10^{-2}$ | $7.43 \cdot 10^{-2}$ | $6.09 \cdot 10^{-2}$ | $5.16 \cdot 10^{-2}$ | $4.85 \cdot 10^{-2}$ |
| 20 | 4.14 | 6.42 | 0.9 | 1.26 | 0.7 | 0.59 | 0.69 |
| 30 | 9.4 | 22.71 | 11.89 | 6.08 | 8.4 | 6.02 | 6.09 |
| 40 | 19.46 | 40.56 | 60.49 | 37.89 | 36.41 | 40.19 | 39.56 |
| Fast $L^1$ | | | | | | | |
| 10 | 0.15 | 0.3 | 0.32 | $8.73 \cdot 10^{-2}$ | 0.54 | 1.33 | 12.2 |
| 20 | 0.59 | 1.23 | 1.72 | 1.53 | 2.9 | 51.09 | 31.87 |
| 30 | 1.57 | 4.87 | 5.27 | 10.28 | 14.63 | 173.4 | 143.23 |
| 40 | 4.69 | 20.1 | 23.47 | 41.28 | 49.31 | 550.65 | 569.87 |

Figure 8: Average run time in seconds for HessUpdate and Fast $L^1$ for regular grids of $\mathbb{R}^2$ of different sizes and different parameters $\alpha$. All solutions converged to within 0.5% of the ground truth value.

Using the same machine, we performed another set of experiments where we compare this time the performance of HessUpdate with the algorithm introduced in [21], which we will refer as Fast $L^1$. Fast $L^1$ is designed to solve the optimal transport problem on a regular grid in $\mathbb{R}^d$ using a finite-volume discretization of mass flow, by applying a primal-dual method for a quadratically-regularized $L^1$ problem. The results of this experiment are in Figure 8

We benchmarked HessUpdate and Fast $L^1$ on regular grids of $\mathbb{R}^2$, with $N \in \{10, 20, 30, 40\}$ nodes per side. This gives a graph $G$ with $N^2 \in \{100, 400, 900, 1600\}$ vertices and $2N(N-1) \in \{180, 760, 1740, 3120\}$ edges. The data $f$ was again constructed by randomly choosing 10% of the vertices, assigning them random values between $-10$ and $10$, and setting the last chosen vertex to contain the appropriate amount of mass to ensure a proper mass balance of sources and sinks. Ten different functions $f$ were generated for each $(N, \alpha)$ pair; average run times are reported.

Our implementation of Fast $L^1$ in MATLAB uses the same parameters as theirs, which are $\mu = \tau = 0.025$ and $\theta = 1$, which guarantees convergence for such values of $N$ according to [21, Theorem 1]. Both primal and dual variables were initialized to 0.

The only difference between the authors' implementation of Fast $L^1$ and ours resides in the stopping condition; to ensure a reasonable comparison between the two algorithms, we first run HessUpdate until convergence or reaching the maximum number of iterations. We then compute the relative error between the objective value obtained by HessUpdate and the ground truth value. Finally, we run Fast $L^1$ until it reaches the same relative error compared to the ground truth value. Both HessUpdate and Fast $L^1$ were limited to single threaded runs, and every attempt was made to optimize the MATLAB implementation of Fast $L^1$ to the point where the "shrink" operation dominates run time.

HessUpdate outperforms Fast $L^1$ for $\alpha > 10^{-2}$ over all grid sizes, and for $\alpha = 10^{-2}$ for small grid sizes (less than 20 nodes per side). For smaller temperatures, Fast $L^1$ is faster for regular grids of $\mathbb{R}^2$. In general, the experiment reveals that Fast $L^1$ may be slightly preferable in the low-temperature regime for grids on $\mathbb{R}^2$, in its specific target case of grid graphs.

# 7 Discussion and Conclusion

Our algorithm for regularized graph transport finds immediate application in software for network flow and related problems. Regularization allows the output for these problems to be *predictable*, since strict convexity guarantees uniqueness of the minimizer. Our technique is efficient and built upon general-purpose tools for sparse linear algebra, which will likely improve in the future.

Several avenues for future research will address remaining theoretical and practical challenges suggested by our work. On the theoretical side, an analog of our analysis may apply to the Beckmann model of transport over $\mathbb{R}^n$ [1, 32], showing how quadratic regularization affects flows in the continuum. In computer science theory, time-complexity analysis related to that in [36, 10] could reveal worst-case behavior of our algorithm in terms of computational operations; note that our Newton step may lead to convergence (with infinite-precision arithmetic) in a finite number of steps. A better characterization of the evolution of the active sets with increasing regularizers and a clean proof of the monotonicity conjecture 1 also would complete our understanding of the regularization. From a practical perspective, future implementations could consider extremely large graphs that cannot fit into memory and/or architectures that require parallel processing with minimal synchronization.

### Acknowledgments

The authors thank Alfred Galichon for introducing and discussing the quadratic regularization, Aleksander Mądry for discussing and helping refine the theoretical characterizations in this paper, as well as the two anonymous referees for their interesting comments that improved theoretical results, and clarity of the paper. J. Solomon acknowledges support of Army Research Office grant W911NF-12-R-0011 ("Smooth Modeling of Flows on Graphs"), M. Essid acknowledges support from NSF grant DMS-1311833. M. Essid is particularly grateful to Alfred Galichon and Robert V. Kohn for their constant help and support discussing theoretical and numerical aspects, as well as Flavien Léger for fruitful discussions.